\documentclass[a4paper,11pt]{article}
\title{Equilateral Sets in Banach Spaces of the form C(K)}
\author{S.K.Mercourakis and G.Vassiliadis}
\date{}
\usepackage{amsmath}
\usepackage{amsfonts}
\usepackage{amssymb}
\usepackage{amscd}
\usepackage{amsthm}
\usepackage{makeidx}
\usepackage{euscript}
\theoremstyle{plain}
\newtheorem{theo}{Theorem}
\newtheorem{lemm}{Lemma}
\newtheorem{prop}{Proposition}
\newtheorem{cor}{Corollary}
\theoremstyle{definition}
\newtheorem{defn}{Definition}
\newtheorem*{claim}{\underline{Claim}}
\newcommand{\beg}{\begin{proof}[Proof]}
\newcommand{\eq}{equilateral $\,$}

\newcommand{\ban}{Banach space $\,$}

\begin{document}
\maketitle

\begin{abstract}
We show that for "most" compact non metrizable spaces, the unit ball of
the \ban $C(K)$ contains an uncountable 2-\eq set. We also give
examples of compact non metrizable spaces $K$ such that the minimum 
cardinality of a maximal \eq set in C(K) is countable.
\end{abstract}


\section*{Introduction}
A subset $S$ of a metric space $(M,d)$ is said to be \eq if there is
a constant $\lambda>0$ such that
$d(x,y)=\lambda$, for $x,y \in S,x \neq y$; we also call such a set 
a $\lambda$-\eq set. An \eq set $S \subseteq M$ is said to be maximal 
if there is no \eq set $B \subseteq M$ with $A \subsetneqq B$.

Equilateral sets have been studied mainly in finite dimensional spaces, see
\cite{P},\cite{SV1} and \cite{S} for a survey on \eq sets. More recently
there are also results on infinite dimensions, see \cite{MV}, \cite{FOSS}
and also on maximal \eq sets, see \cite{SV2}.

In this paper we study \eq sets in Banach spaces of the form C(K), where $K$ is a compact
space. The paper is divided into two sections. In the first section we introduce 
the combinatorial concept of a linked family of pairs of a set $\Gamma$; using
this concept we characterize those compact spaces $K$, such that the unit ball
of $C(K)$ contains a $(1+\varepsilon)$-separated (equivalently: a 2-equilateral) set
of a given cardinality (Theorem 1). Then we show that in "most" cases a compact
non metrizable space $K$ admits an uncountable linked family of closed pairs and 
hence its unit ball contains an uncountable 2-\eq set (Theorem 2).

In the second section we focus on maximal \eq sets on the space
$C(K)$. Following \cite{SV2} (Definition 2), given a normed space E,
we denote by $m(E)$ the minimum cardinality of a maximal \eq set in E.
The main results here are the following: for every infinite locally compact
space $K$ we have $m(C_0(K)) \ge \omega$ (Theorem 3) (thus in particular, 
$m(C(K))=\omega$, for any infinite compact metric space $K$). For every infinite
product $K=\prod\limits_{\gamma \in \Gamma}^{} K_\gamma$ of nontrivial compact
metric spaces, $m(C(K))=|\Gamma|$ (Theorem 4). Then using proper linked families
of pairs of $\Bbb{N}$, we give a variety of examples of compact 
non metrizable spaces $K$ (including scattered compact and the Stone-$\check{C}$ech
compactification $\beta \Gamma$ of any infinite discrete set $\Gamma$) such that
$m(C(K))=\omega$ (Theorems 6, 7 and Corol.5).

If E is any (real) \ban then $B_X$ denotes its closed unit
ball. If $K$ is any compact Hausdorff
space, then $C(K)$ is the \ban of all continuous real functions on
$K$ endowed with the supremum norm $\|\cdot\|_\infty$.

\section*{Linked families and equilateral sets in Banach spaces of the form $C(K)$}

In this section we introduce the concept of a linked family of pairs
of a set $\Gamma$ and then use this, in order to investigate the existence 
of \eq sets in $C(K)$, where $K$ is any compact (non metrizable) space.

\begin{defn}
Let $\mathcal{F}=\{(A_\alpha,B_\alpha): \alpha \in \mathcal{A}\}$
be a family of pairs of a nonempty set $\Gamma$. We say that this family 
is \underline{linked} (or \underline{intersecting}) if
\begin{enumerate}
\item $A_\alpha \cap B_\alpha=\emptyset$, for $\alpha \in \mathcal{A}$
\item $A_\alpha \cup B_\alpha \neq \emptyset$, for $\alpha \in \mathcal{A}$ and
\item for $\alpha,\beta \in \mathcal{A}$ with $\alpha \neq \beta$ we
have, either $A_\alpha \cap B_\beta \neq \emptyset$ or $A_\beta
\cap B_\alpha \neq \emptyset$.
\end{enumerate}
\end{defn}

If we replace condition (2) with the stronger one:  
($2^{'}$) $A_\alpha \neq \emptyset \neq B_\alpha$, for $\alpha \in \mathcal{A}$
we shall say that $\mathcal{F}$ is a linked family of nonempty pairs.

We note the following easily verified facts:
\begin{enumerate}
\item If $\alpha \neq \beta \in \mathcal{A}$ then $A_\alpha \neq A_\beta$
and $B_\alpha \neq B_\beta$ (and hence)
\item there is at most one $\alpha \in \mathcal{A}$ such that $A_\alpha=\emptyset$
and at most one $\beta \in \mathcal{A}$ such that $B_\beta=\emptyset$.
\item If $\mathcal{F}$ is a linked family of nonempty pairs of the set $\Gamma$,
then the family $\mathcal{F} \cup \{(\Gamma,\emptyset),(\emptyset,\Gamma)\}$
is a linked family of pairs of $\Gamma$.
\end{enumerate}

\noindent \textbf{Examples 1} 
\begin{enumerate}
\item Let $\{A_\alpha:\alpha \in \mathcal{A}\}$ be a family of distinct 
subsets of the set $\Gamma$. Then the family
$\{(A_\alpha,\Gamma \setminus A_\alpha):\alpha \in \mathcal{A}\}$ is
linked. Assuming furthermore that $\emptyset \neq A_\alpha \neq \Gamma$, for
$\alpha \in \mathcal{A}$, we get that $\mathcal{F}$ is a linked family of nonempty pairs.
It follows in particular that the family $\{(A,\Gamma \setminus A):A \subseteq \Gamma \}$
(resp. the family $\{(A,\Gamma \setminus A):\emptyset \neq A \subsetneqq \Gamma \}$ ) 
is linked (resp. a linked family of nonempty pairs).
\item Let $K$ be a compact totally disconnected space, then the family
$\{(V,K \setminus V):V$ is a clopen subset of $K \}$ is linked.
\item The family $\mathcal{F}=\{(\{1\},\{2\}),(\{2\},\{3\}),(\{3\},\{1\})\}$ is a linked
family of nonempty pairs of the set $\Gamma=\{1,2,3\}$.
\end{enumerate}

In section 2 we shall present much more examples of linked families. Now we are going to
examine the interrelation between the concepts of linked families and \eq sets
in Banach spaces of the form $C(K)$, where $K$ is a compact Hausdorff space.

\begin{lemm} Let $K$ be a compact Hausdorff space and $S \subseteq [0,1]^K \cap C(K)$.
Set $A_f=f^{-1}(\{0\})$ and $B_f=f^{-1}(\{1\})$, for $f \in S$. Then the following
are equivalent:\\
(a) The family $\mathcal{F}=\{(A_f,B_f): f \in S \}$ of (closed) pairs of $K$  is linked.\\
(b) The set $S$ is 1-\eq in $C(K)$.
\end{lemm}

\beg (a) $\Rightarrow$ (b) Let $f,g \in S$ with $f \neq g$; clearly
$0<||f-g||_\infty \le 1$. Since we either have $A_f \cap B_g \neq \emptyset$
or $A_g \cap B_f \neq \emptyset$, there is a $t_0 \in K$ such that
$|f(t_0)-g(t_0)|=1$, hence $||f-g||_\infty=1$.\\
(b) $\Rightarrow$ (a) If $f \neq g \in S$, then $||f-g||_\infty=1$;
so by the compactness of $K$, there is a $t_0 \in K$ such that 
$|f(t_0)-g(t_0)|=||f-g||_\infty=1$. Since $0 \le f(t_0),g(t_0) \le 1$
we get that $\{f(t_0),g(t_0)\}=\{0,1\}$. Therefore, either $t_0 \in
A_f \cap B_g$ or $t_0 \in
A_g \cap B_f$ and $\mathcal{F}$ is as required.
\end{proof}

\noindent \textbf{Note:} Since there is at most one $f_0 \in S$
with $A_{f_0}=\emptyset$ ($\Leftrightarrow inf(f_0)>0$) and at most 
one $g_0 \in S$ with $B_{g_0}=\emptyset$ ($\Leftrightarrow ||g_0||_\infty <1$),
we get that the family $\{(A_f,B_f): f \in S \setminus \{f_0,g_0\} \}$
is a linked family of nonempty closed pairs of $K$ and the set $S \setminus \{g_0\}$
is a subset of the positive part $S^{+}_{C(K)}$ of the unit sphere $S_{C(K)}$ of the
space $C(K)$.

\begin{lemm}
Let $\mathcal{F}=\{(A_\alpha,B_\alpha): \alpha \in \mathcal{A}\}$ be a linked family
of closed pairs of the compact space $K$. Then we can associate with $\mathcal{F}$
a 1-\eq subset $S$ of $C(K)$ with $|S|=|\mathcal{A}|$ and $S \subseteq [0,1]^K \cap C(K)$.
\end{lemm}

\beg Let $\alpha \in \mathcal{A}$; we distinguish the following cases for the pair
$(A_\alpha,B_\alpha)$:\\
(I) $A_\alpha \neq \emptyset \neq B_\alpha$. We consider a Urysohn function
$f_\alpha:K \rightarrow [0,1]$ so that $f_\alpha(x)=0$ for $x \in A_\alpha$ and
$f_\alpha(x)=1$ for $x \in B_\alpha$; clearly $inf(f_\alpha)=0<||f_\alpha||_\infty=1$.\\
(II) Assume that $A_\alpha=\emptyset$, thus $B_\alpha \neq \emptyset$. If 
$B_\alpha \neq K$, pick $t_0 \in K \setminus B_\alpha$ and consider a Urysohn
function $f_\alpha:K \rightarrow [0,1]$ so that $f_\alpha \diagup B_\alpha=1$
and $f_\alpha(t_0)=0$. In case when $B_\alpha=K$, we let $f_\alpha=1$ on $K$.\\
(III)  Assume that $B_\alpha=\emptyset$, thus $A_\alpha \neq \emptyset$.
This case is similar to case (II). So we consider a Urysohn function
$f_\alpha:K \rightarrow [0,1]$ so that $f_\alpha \diagup A_\alpha=0$ and
$f_\alpha(t_0)=1$ for some $t_0 \in K \setminus A_\alpha$, if $A_\alpha \neq K$
and define $f_\alpha$ to be the constant zero function in case when $A_\alpha=K$.

Now set $A_\alpha^{'}=f_\alpha^{-1}(\{0\})$ and  $B_\alpha^{'}=f_\alpha^{-1}(\{1\})$,
for $\alpha  \in \mathcal{A}$. Since $A_\alpha^{'} \cap B_\alpha^{'} = \emptyset$ , 
$A_\alpha \subseteq A_\alpha^{'}$ and $B_\alpha \subseteq B_\alpha^{'}$ 
for $\alpha  \in \mathcal{A}$, we have that the family $\{(A_\alpha^{'},B_\alpha^{'}):
\alpha  \in \mathcal{A}\}$ is a linked family of closed pairs of the space $K$, hence by
Lemma 1 the set $S=\{f_\alpha:\alpha \in \mathcal{A}\}$ is a 1-\eq subset of $[0,1]^K \cap C(K)$.
\end{proof}

\noindent \textbf{Remarks 1} (1) If in the proof of Lemma 2 we consider
(as we may) continuous functions $f_\alpha:K \rightarrow [-1,1]$ such that
$f_\alpha \diagup A_\alpha=1$ and $f_\alpha \diagup B_\alpha=-1$, then the
set $\{f_\alpha:\alpha \in \mathcal{A}\}$ is a 2-\eq subset of the unit ball of $C(K)$.\\
(2) Let $\mathcal{F}=\{(A_\alpha,B_\alpha):\alpha \in \mathcal{A}\}$ be a family
of disjoint pairs of a set $\Gamma$. Set $\bar{\mathcal{F}}=
\{(\bar{A}_\alpha,\bar{B}_\alpha):\alpha \in \mathcal{A}\}$ where 
$\bar{A}_\alpha=cl_{\beta \Gamma}A_\alpha$, $\bar{B}_\alpha=cl_{\beta \Gamma}B_\alpha$
and $\beta \Gamma$ is the Stone-$\check{C}$ech compactification  of the discrete set $\Gamma$.
Then it is easy to see that $\mathcal{F}$ is a linked family of (nonempty) pairs of $\Gamma$
iff $\bar{\mathcal{F}}$ is a linked family of (nonempty) pairs of $\beta \Gamma$.\\
(3) Let $K$ be a compact space and $S \subseteq [0,1]^K \cap C(K)$ be a 1-\eq set. We consider
the linked family $\mathcal{F}=\{(A_f,B_f):f \in S\}$ given by Lemma 1. Then it is 
not difficult to verify that $\mathcal{F}$ is a maximal linked family of closed pairs of $K$
iff the set $S$ is a maximal (with respect to inclusion) 1-\eq subset of $[0,1]^K \cap C(K)$,
endowed with the norm metric (the equilateral set $S$ is not necessarily maximal in the space
$C(K)$, see Remark 5(3) ).

\begin{theo}
Let $K$ be a compact Hausdorff space and $\alpha$ be an infinite
cardinal. The following are equivalent:
\begin{enumerate}
\item The unit ball $B_{C(K)}$ of $C(K)$ contains a $\lambda$-\eq set with $\lambda>1$,
of size $\alpha$.
\item The unit sphere $S_{C(K)}$ (resp. the positive part of the unit sphere $S^{+}_{C(K)}$)
of $C(K)$ admits a 2-\eq (resp. a 1-equilateral) set of size $\alpha$.
\item The unit ball $B_{C(K)}$ of $C(K)$ contains a $(1+\varepsilon)$-separated set, for some $\varepsilon>0$,
of size equal to $\alpha$.
\item There exists a linked family of closed (nonempty) pairs in $K$ of size
equal to $\alpha$.
\end{enumerate}
\end{theo}

\beg $(2) \Rightarrow (1)$ Let $S$ be a 1-\eq subset of $S^{+}_{C(K)}$ with $|S|=\alpha$.
Then by Lemma 1, $\mathcal{F}=\{(A_f,B_f):f \in S\}$ is a linked family of closed pairs
of $K$ with $|\mathcal{F}|=\alpha$. Therefore, by Lemma 2 and Remark 1(1) $\mathcal{F}$
defines a 2-\eq set contained in the unit ball of $C(K)$.

$(1) \Rightarrow (3)$ is obvious.

$(3) \Rightarrow (4)$ Let $D \subseteq B_{C(K)}$ be a $(1+\varepsilon)$-separated set
($\varepsilon>0$), with $|D|=\alpha$. We may assume that $||f||_\infty=1$, for
$f \in D$. We define $A_f=f^{-1}([-1,-\frac{\varepsilon}{2}])$ and 
$B_f=f^{-1}([\frac{\varepsilon}{2},1])$, for $f \in D$; clearly $A_f \cup B_f \neq \emptyset$.
Let $f,g \in D$ with $f \neq g$, so there is a $t_0 \in K$ such that $||f-g||_\infty=
|f(t_0)-g(t_0)| \ge 1+\varepsilon$. Assume without loss of generality that $f(t_0)<g(t_0)$;
then we have, $f(t_0) \le -\frac{\varepsilon}{2}$ and $g(t_0) \ge \frac{\varepsilon}{2}$, that is,
$A_f \cap B_g \neq \emptyset$. For, suppose otherwise, then we would either have $f(t_0) > -\frac{\varepsilon}{2}$
or $g(t_0) <\frac{\varepsilon}{2}$. Assuming that $f(t_0) > -\frac{\varepsilon}{2}$ we get that 
$-\frac{\varepsilon}{2}<f(t_0)<g(t_0) \le 1$, hence $g(t_0)-f(t_0)<1+ \frac{\varepsilon}{2}$, a contradiction.

In a similar way we get a contradiction assuming that $g(t_0) <\frac{\varepsilon}{2}$. It follows that the
family $\mathcal{F}=\{(A_f,B_f):f \in D\}$ defined above is a linked family of closed pairs in $K$ of size
$\alpha$.

$(4) \Rightarrow (2)$ This implication is a direct consequence of Lemma 2.\\
The proof of the Theorem is complete.
\end{proof}

Let $K$ be an infinite compact space; as is well known the \ban $c_0$ is isometrically
embeded in $C(K)$, hence the assertions of Theorem 1 hold true for $\alpha=\omega$.
The following questions are open for us:

\underline{Questions.} Let $K$ be a compact Hausdorff non metrizable space.
\begin{enumerate}
\item Does there exist an uncountable $(1+\varepsilon)$-separated $D \subseteq B_{C(K)}$?
Does there exist (at least) an uncountable $D \subseteq B_{C(K)}$ such that $f \neq g \in D 
\Rightarrow ||f-g||_\infty>1$?

Note that the unit ball of every infinite dimensional \ban contains an infinite $(1+\varepsilon)$-separated
set, see \cite{EO}.
\item Does the space $C(K)$ contain an uncountable \eq set?
\end{enumerate} 

We note that, regarding question (1), by transfinite induction it can be shown that there is
an uncountable $D \subseteq S^{+}_{C(K)}$ such that $f \neq g \in D 
\Rightarrow ||f-g||_\infty \ge 1$.

However we can show that in "most" cases the answer to the above questions
is positive. For this purpose
we recall that a (Hausdorff and completely regular) topological
space $X$ is said to be:

(i) \underline{hereditarily Lindel\"{o}f} (HL) if every subspace $Y$
of $X$ is Lindel\"{o}f. It is well known that a space $X$ is HL iff
there is no uncountable right separated family in $X$; that is, a
family $\{t_\alpha:\alpha<\omega_1\} \subseteq X$ such that
$t_\alpha \notin cl_X \{t_\beta:\alpha<\beta<\omega_1\}$ for
$\alpha<\omega_1$ and

(ii) \underline{hereditarily separable} (HS) if every subspace $Y$
of $X$ is separable. It is also well known that a space $X$ is HS iff
there is no uncountable left separated family in $X$; that is, a
family $\{t_\alpha:\alpha<\omega_1\} \subseteq X$ such that
$t_\alpha \notin cl_X \{t_\beta:\beta<\alpha \}$ for $1 \le
\alpha<\omega_1$ (see \cite{HMVZ} p. 151).

We are going to use the following standard\\
\textbf{Fact.} A compact space $K$ is HL if and only if it is
perfectly normal (i.e. each closed subset of $K$ is $G_\delta$).

\begin{theo}
Let $K$ be a compact space. If $K$ satisfies one of the following
conditions, then $K$ admits an uncountable linked family of closed
pairs (and hence by Theorem 1 the unit ball of $C(K)$ contains an uncountable 2-\eq set).
\begin{enumerate}
\item There exists a closed subset $\Omega$ of $K$ admitting
uncountably many relatively clopen sets (in particular $\Omega$
is non metrizable and totally disconnected).
\item $K$ is non hereditarily Lindel\"{o}f.
\item $K$ is non hereditarily separable.
\item $|K|>c=2^\omega$(=the cardinality of continuum).
\item $K$ admits a Radon probability measure of uncountable type.
\end{enumerate}
\end{theo}

\beg

(1) Let $\mathcal{B}$ be any uncountable family of clopen sets
in $\Omega \subseteq K$. Then clearly the family  
$\mathcal{F}=\{(V,\Omega \setminus V):V \in \mathcal{B} \}$ is an uncountable linked family of closed
pairs in $K$. (It is clear that condition (1) can be stated as follows: there is a closed subset 
$\Omega \subseteq K$ such that the unit ball of $C(\Omega)$ has uncountably many extreme points).

(2) Let $\{t_\alpha:\alpha<\omega_1\} \subseteq K$ be an uncountable
right separated family. Set $A_\alpha=\{t_\alpha\}$ and  $B_\alpha=cl_K \{t_\beta:\alpha<\beta<\omega_1\}$
for $\alpha<\omega_1$. Then it is easy to see that the family
$\{(A_\alpha,B_\alpha):\alpha<\omega_1\}$ is a linked family of closed (nonempty) pairs of $K$.

(3) Since $K$ is non HS, there exists an uncountable left separated
family in $K$ and the proof is similar to that of the previous case.

(4) This follows from (2), since if $|K|>c$ then $K$ is not HL. Indeed, any compact HL space is
first countable (each point set of $K$ is $G_\delta$ by the Fact
preceding the Theorem). By a classical result of Archangel'skii
each compact first countable space has cardinality $\le c$.

(5) Let  $\mu \in P(K)$ be a Radon probability measure on $K$ 
of uncountabe type (i.e. $dimL_1(\mu) \ge \omega_1$).
We consider (as we may) an uncountable stochastically independent
family $\{\Gamma_\alpha:a<\omega_1\}$ of $\mu$-measurable subsets
of $K$. This means that \[\mu(\bigcap\limits_{k=1}^n \varepsilon_k \Gamma_{\alpha_k})=\frac{1}{2^n}, \,\,
\alpha_1<\cdots<\alpha_n, \;\;\textrm{and} \,\, \varepsilon_1, \cdots ,\varepsilon_n \in \{-1,1\}\]
where, if $A \subseteq K$ we let $1 \cdot A=A$ and $(-1) \cdot A=K \setminus A$. By the regularity
of the measure $\mu$ we can find compact sets $A_\alpha \subseteq \Gamma_\alpha$ and $B_\alpha \subseteq
K \setminus \Gamma_\alpha$, for each $\alpha<\omega_1$ such that $\mu(\Gamma_\alpha \setminus A_\alpha)<\frac{1}{8}$
and $\mu((K \setminus \Gamma_\alpha) \setminus B_\alpha)<\frac{1}{8}$ \,\, (1).

\begin{claim}
The family of closed pairs $\{(A_\alpha,B_\alpha):\alpha<\omega_1\}$ is linked. 
\end{claim}

\underline{Proof of the Claim:} Let $\alpha<\beta<\omega_1$; we are going to show the 
stronger property $\mu(A_\alpha \cap B_\beta)>0$ and $\mu(A_\beta \cap B_\alpha)>0$.
Assume that $\mu(A_\alpha \cap B_\beta)=0$; it then follows from (1) that 
$\mu(\Gamma_\alpha \setminus A_\alpha)=\mu(\Gamma_\alpha)-\mu(A_\alpha)=
\frac{1}{2}-\mu(A_\alpha)<\frac{1}{8}$, thus $\mu(A_\alpha)>\frac{1}{2}-\frac{1}{8}=\frac{3}{8}$.
Also $\mu((K \setminus \Gamma_\beta) \setminus B_\beta)=\mu(K \setminus \Gamma_\beta)-\mu(B_\beta)=
\frac{1}{2}-\mu(B_\beta)<\frac{1}{8}$, thus $\mu(B_\beta)>\frac{1}{2}-\frac{1}{8}=\frac{3}{8}$.

Therefore $\mu(\Gamma_\alpha \cup (K \setminus \Gamma_\beta)) \ge \mu(A_\alpha \cup B_\beta)
=\mu(A_\alpha)+\mu(B_\beta)>\frac{3}{8}+\frac{3}{8}=\frac{3}{4}$, a contradiction because
 $\mu(\Gamma_\alpha \cup (K \setminus \Gamma_\beta))=\frac{3}{4}$.

In a similar way we get that $\mu(A_\beta \cap B_\alpha)>0$ and the proof of the Claim is complete.

\end{proof}

The above Theorem has some interesting consequences. If $K$ is any compact space, then $P(K)$
denotes the set of Radon probability measures on $K$. Recall that both spaces $P(K)$ and $B_{C(K)^*}$
are compact with the weak$^*$ topology.

\begin{cor}
Let $K$ be any compact non metrizable space. Denote by
$\Omega$ any of the compact spaces $K \times K,P(K)$ and
$B_{C(K)^*}$. Then the unit ball of $C(\Omega)$ contains an uncountable 2-\eq
set.
\end{cor}

\beg The compact space $\Omega$ is not HL. Indeed, if $\Omega=K
\times K$ , then since $K$ is not metrizable, its diagonal
$\Delta=\{(x,x):x \in K\}$ is closed but not $G_\delta$ subset of
$\Omega$ (by a classical result if the diagonal of a compact space
$K$ is a $G_\delta$ subset of $K \times K$ then $K$ is metrizable). So
by the Fact before Th.1 the space $K \times K$ is not HL.

Let $\Omega=P(K)$. We consider the continuous map $\Phi:K \times K
\rightarrow P(K):\Phi(x,y)=\frac{1}{2} \delta_x+\frac{1}{2}
\delta_y$ ($\delta_x$ is the Dirac measure at $x \in K$). Then
$\Delta=\Phi^{-1}(\{\delta_x:x \in K\})$. If the space $P(K)$ were HL, then
by the Fact above $K$ would be a closed $G_\delta$ subset of $P(K)$,
therefore $\Delta$ would be a $G_\delta$ subset of $K \times K$, a
contradiction.

If $\Omega=B_{C(K)^*}$ then since $P(K)$ is a weak$^*$ closed
subset of $\Omega$ we get that $\Omega$ is not HL.
\end{proof}

\begin{cor}
Let $X$ be a nonseparable Banach space. Denote by $\Omega$ its
closed dual unit ball $B_{X^*}$ with the weak$^*$ topology.
Then the unit ball of $C(\Omega)$ contains an uncountable 2-\eq set.
\end{cor}

\beg Using transfinite induction and Hahn-Banach Theorem we may
construct for each $\varepsilon>0$ two long sequences
$\{x_\alpha:\alpha<\omega_1\} \subseteq B_X$ and
$\{f_\alpha:\alpha<\omega_1\} \subseteq (1+\varepsilon)B_{X^*}$
satisfying $f_\beta(x_\alpha)=0$ for $\beta>\alpha$ and
$f_\alpha(x_\alpha)=1$ for $\alpha<\omega_1$ (see the Fact 4.27 of
\cite{HMVZ}). It is easy to see that the sequence
$\{f_\alpha:\alpha<\omega_1\}$ is right separated in the compact
space $(1+\varepsilon)B_{X^*}$. So the space
$(1+\varepsilon)B_{X^*}$ is not HL and the same is valid for
$\Omega=B_{X^*}$.
\end{proof}

In the sequel we relate the concept of a linked family of closed
pairs with the known concept of a weakly separated subspace of
some topological space (\cite{BGT},\cite{LAT}).

A subspace $Y$ of a topological space $X$ is said to be
\underline{weakly separated} if there are open sets $U_y,y \in Y$
in $X$ such that $y \in U_y \,\, \forall y \in Y$ and whenever
$y_1 \neq y_2 \in Y$ we either have $y_1 \notin U_{y_2}$ or $y_2 \notin U_{y_1}$.\\
We note the following easily verified facts:

(i) If $Y=\{t_\alpha:\alpha<\omega_1\}$ is any right (resp. left)
separated family in the topological space $X$, then $Y$ is an
uncountable weakly separated subspace of $X$; we say in this case
that $Y$ is an uncountable right (resp. left) separated subspace of
$X$.

(ii) Let $Y$ be any weakly separated subspace of $X$ by the family of
open sets $U_y, y \in Y$. Then the family $\{(\{y\},X
\setminus U_y):y \in Y\}$ is a linked family of closed pairs in $X$.

As we shall see, linked families of closed pairs in a topological
space $X$ can be interpreted as a special kind of weakly separated
subspaces in $expX$, the hyperspace of  closed nonempty subsets of
$X$ endowed with the Vietoris topology. If $G_1, \dots,G_n$ are
subsets of $X$ we denote $<G_1,\dots,G_n>=\{F \in expX:F \subseteq
\bigcup\limits_{k=1}^n G_k \,\, \textrm{and} \,\, F \cap G_k \neq \emptyset \,\, \forall
k=1,2,\dots,n \}$. The Vietoris topology on $expX$ has as base the sets of
the form $<G_1,\dots,G_n>$, where $G_1,\dots,G_n$ are open subsets
of $X$.

We shall say that a weakly separated subspace $Y$ of $expX$ is
separated by \underline{open subsets of $X$}, if the sets $U_y, y
\in Y$ of the definition above are of the form $U_y=<V_y>,y \in
Y$, where $V_y$ are open subsets of $X$. This is equivalent to both:
$y \subseteq V_y$ for $y \in Y$ and if $y_1 \neq y_2 \in Y$, then either 
$y_1 \nsubseteq V_{y_2}$ or $y_2 \nsubseteq
V_{y_1}$. 

More exactly we have the following (easy) Proposition, the proof of which 
is left to the reader.

\begin{prop}
Let $X$ be a topological space and $\kappa$ any cardinal. The
following are equivalent:
\begin{enumerate}
\item $X$ admits a linked family of closed pairs of cardinality
$\kappa$.
\item $expX$ admits a weakly separated subspace by open subsets of
$X$ of cardinality $\kappa$.
\end{enumerate}
\end{prop}

\noindent \textbf{Remarks 2} (1) As was shown by Todorcevic
assuming Martin's Axiom and the negation of the continuum
hypothesis, if $K$ is compact and non metrizable then the space
$C(K)$ admits an uncountable (bounded) biorthogonal system
(\cite{T}, Th.11). So by using Theorem 3 of \cite{MV}, the space $C(K)$ can be
given an equivalent norm that admits an uncountable \eq set.

(2) It is consistent with ZFC to assume that there exists a compact
non metrizable space $K$ having no uncountable weakly separated
subspace (see \cite{BGT}). The space $K$ constructed there, among
its many interesting properties, is totally disconnected and hence admits  an
uncountable linked family of closed (and open) pairs.

(3) Let $K$ be a compact non metrizable space. Then the hyperspace
$expK$ of $K$ is not HL. Actually its closed subspace $[K]^{\le
2}=\{A \subseteq K:|A| \le 2\}$ is not HL. (The proof is similar
to the proof that the space $(P(K),w^*)$ is not HL; we consider
the continuous map $\Phi:K \times K \rightarrow
expK:\Phi(x,y)=\{x,y\}$ and note that $\Phi(K \times K)=[K]^{\le
2}$). It follows that there is an uncountable right separated
subspace $Y=\{F_\alpha:\alpha<\omega_1\}$ of $expK$, but it is not
clear whether $Y$ is (or another uncountable weakly separated
subspace of $expK$ can be chosen so as to be) separated by open
subsets of $K$.

\section*{Maximal equilateral sets in Banach spaces of the form $C(K)$}

Our goal here is the study of maximal \eq sets of minimum cardinality,
mainly in Banach spaces of the form $C(K)$. As we shall see, proper 
linked families of pairs of $\Bbb{N}$ play a key role.

\begin{defn}
Let $(M,d)$ be a metric space. We define, for $x \in M$, $m(M,x)=min\{|A|: x \in A$
and $A$ is a maximal \eq set in M$\}$. We also define $m(M)=min\{|A|:$ $A$ is a maximal
\eq set in M$\}$
\end{defn}

It is clear that $m(M)=min\{m(M,x): x \in M\}$.

\begin{lemm}
Let $(X,|| \cdot ||)$ be a normed space. Then we have $m(X)=m(S_X \cup \{0\},0)$.
\end{lemm}

\beg Let $A \subseteq X$ be any maximal \eq set in $X$. Assume that $A$ is $\lambda$-equilateral.
Let $x_0 \in A$; then the set $B=\{\frac{1}{\lambda} (x-x_0):x \in A\}$ is a 1-\eq 
subset of $S_X \cup \{0\}$ containing the point 0. Note that $|B|=|A|$ and that B is a
maximal \eq set (in $X$ and hence) in the metric space  $S_X \cup \{0\}$.

In the converse direction, consider any maximal \eq subset B of the metric space  $S_X \cup \{0\}$
with $0 \in B$. Then clearly B is 1-equilateral. We claim that B is a maximal \eq subset of $X$; indeed, if $x \in X$
with $x \notin B$ such that $B \cup \{x\}$ is \eq then $1=||x-0||=||x||$, so $x \in S_X$ which
contradicts the maximality of B in the metric space $S_X \cup \{0\}$.
\end{proof}

\begin{lemm}
Let $(X,|| \cdot ||)$ be a normed space. Then we have $m(B_X) \le m(X) \, (=m(S_X \cup \{0\},0))$.
\end{lemm}

\beg By the (method of proof of) the previous Lemma any maximal \eq set in $X$ gives
rise to a maximal \eq set in $X$ of the same cardinality, contained in $S_X \cup \{0\} \subseteq B_X$,
so we are done.
\end{proof}

\noindent \textbf{Remarks 3} (1) Swanepoel and Villa have shown in \cite{SV2} the following result, 
generalizing an example of Petty \cite{P}:\\
If $X$ is any \ban with $dimX \ge 2$ having a norm which is G\^{a}teaux differentiable at some point
and $Y=(X \oplus \Bbb{R})_1$, then we have $m(Y)=4$.\\ 
(Their proof is based on the following simple result:
Let $X$ be any normed space with $dimX \ge 2$ and also let $x,u \in S_X$ such that $||u-x||=||u+x||=2$. Then the unit ball
of the subspace $Z=<u,x>$ of $X$ is the parallelogram with vertices $\pm u,\pm x$.) One can easily check  that the result 
of Swanepoel and Villa can be generalized (by the method of its proof) as follows:\\
If $dimX \ge 2$ and the norm of $X$ is either strictly convex or G\^{a}teaux differentiable at some point, then we have
$m(Y)=4$ and $m(B_X)=2(=m(S_X))$.\\
(2) For the Hilbert space $X=\ell_2$ clearly we have $m(X)=\omega$ and since the norm of $X$ is strictly convex, we get
from the preceding remark that $m(B_X)=2$. So the inequality in Lemma 4 can be strict.\\
(3) Let $\Gamma$ be any set with $|\Gamma| \ge 2$ and let $|| \cdot||$ be an equivalent strictly convex norm
on the \ban $\ell_1(\Gamma)$, see \cite{DGZ}. Now set $X=(\ell_1(\Gamma),|| \cdot ||)$, then we get from Remark 3(1) that
$m(X \oplus \Bbb{R})_1=4$.

Now we are going to generalize the following result of Swanepoel and Villa in \cite{SV2}: If $d \in \Bbb{N}$
then $m(\ell_\infty^d)=d+1$.\\
Let $X$ be a locally compact Hausdorff space and $C_0(X)$ be the \ban (endowed with supremum norm) of all 
continuous functions $f:X \rightarrow \Bbb{R}$ with the property that, $\forall \varepsilon>0 \,\, \exists K \subseteq X$
compact: $|f(x)|< \varepsilon$, for all $x \in X \setminus K$. As we know $C_0(X)$ is the completion of the space
$C_c(X)$ of continuous functions $f:X \rightarrow \Bbb{R}$ with compact support.
We note the following facts:\\
(1) If $f,g \in C_0(X)$, then $max\{f,g\}$ and $min\{f,g\}$ belong to $C_0(X)$.\\
(2) If $A$ is any finite nonempty subset of $X$ and $g:A \rightarrow \Bbb{R}$ (resp. $g:A \rightarrow [0,1]$)
is any function, then there is a continuous extension $f:X \rightarrow \Bbb{R}$ (resp. $f:X \rightarrow [0,1]$)
of $g$, which has compact support (the proof of this fact uses Urysohn's Lemma).

\begin{theo}
Let $X$ be any infinite locally compact Hausdorff space. Then we have $m(C_0(X)) \ge \omega$.
\end{theo}

\beg We shall show that each finite \eq subset of $C_0(X)$ can be extended. So let $S=\{f_1, \ldots, f_n \}, \,
n \ge 2$ be any 1-\eq set in $C_0(X)$. Since $|f_k(x)-f_l(x)| \le 1$ for all $k \neq l \le n$ and $x \in X$, we
may assume that $S \subseteq [0,1]^X \cap C_0(X)$. (Indeed, set $f(x)=min\{f_k(x):1 \le k \le n\}$ for $x \in X$,
then the function $f \in C_0(X)$ and $0 \le f_k(x)-f(x) \le 1$ for $k \le n$ and $x \in X$. So the set 
$\{g_k=f_k-f:1 \le k \le n\}$ is a 1-\eq subset of $[0,1]^X \cap C_0(X)$).

We consider any finite subset $A$ of $X$ with $|A| \ge n$, such that:\\
(1) $\forall k \le n \,\, \exists t_k \in A$ with $|f_k(t_k)|=||f_k||_\infty$ and\\
(2) $\forall k \neq l \le n \,\, \exists t=t(k,l) \in A: ||f_k-f_l||_\infty=|(f_k-f_l)(t)|=1$.

Then the set $\{h_k=f_k \diagup A:1 \le k \le n\}$ is 1-\eq in the space $\ell_\infty(A)$ and
since $n \le |A|< \omega$, it can be extended  on $\ell_\infty(A)$ to a 1-\eq set with at least
$|A|+1 \ge n+1$ elements (see Prop.12 in \cite{SV2}). Let $h \in \ell_\infty(A)$ taking values
in [0,1] such that $||h-h_k||_\infty=1$ for $k \le n$. Then by using Fact (2) mentioned above, we can
find a continuous extension $f:X \rightarrow [0,1]$ of $h$ on $X$ having compact support. It is obvious
that the set $S \cup \{f\}$ is a 1-\eq set in $C_0(X)$, so we are done.
\end{proof}

Let $\Gamma$ be an infinite set endowed with discrete topology. Then $c_0(\Gamma)$ is the space of
all functions $f:\Gamma \rightarrow \Bbb{R}$ that vanish at infinity. We shall show that the number
$m(c_0(\Gamma))$ is as big as possible.

\begin{prop}
Let $\Gamma$ be an infinite set, then $m(c_0(\Gamma))=|\Gamma|$
\end{prop}

\beg We claim that each \eq subset $S$ in $c_0(\Gamma)$ with $|S|<|\Gamma|$ can be extended. If $\Gamma$
is countable, then $S$ is finite and can be extended by the previous theorem. So assume $\Gamma$ is uncountable.
It is also clear by Lemma 3 that we may assume $S \subseteq B_{c_0(\Gamma)}$ and that it is 1-equilateral.
Set $\Delta=\cup \{suppx:x \in S\}$; since $|S|<|\Gamma| \ge \omega_1$ and each element of $c_0(\Gamma)$
has at most countable support, we get that $|\Delta|<|\Gamma|$. Let $\gamma_0 \in \Gamma \setminus \Delta$,
then it is easy to see that the set $S \cup \{e_{\gamma_0}\}$ ($e_{\gamma_0}$ is the $\gamma_0$-member of the 
usual basis of $c_0(\Gamma)$) is 1-equilateral. Now we can proceed by transfinite induction, using the above claim
to show that $m(c_0(\Gamma))=|\Gamma|$. We omit the details of this (easy) proof.
\end{proof}

Let $K$ be any infinite compact metric space, then the \ban $C(K)$ is separable and by Theorem 3 we get that
$m(C(K))=\omega$. This result can be generalized as follows:

\begin{theo}
Let $\{X_\gamma:\gamma \in \Gamma\}$ be an infinite family of compact metric spaces, so that
$|X_\gamma| \ge 2$, for all $\gamma \in \Gamma$. Set $X=\prod\limits_{\gamma \in \Gamma} X_\gamma$, then we have
$m(C(K))=|\Gamma|$.
\end{theo}

The main tool for proving Theorem 4 is the following

\begin{lemm}
Let $X,Y$ be compact spaces such that $X$ is metrizable and $|Y| \ge 2$. Let $S$ be any \eq set
in $C(X)$. Then there is a linear isometry extension operator $T:C(X) \rightarrow C(X \times Y)$,
so that the set $T(S)$ can be extended to an \eq set in $C(X \times Y)$.
\end{lemm}

\beg Assume without loss of generality that $S$ is a 1-\eq set such that $S \subseteq S_{C(X)} \cup \{0\}$.
Consider two distinct elements, say $a, b$ of the space $Y$ and set $E=C(X \times \{a,b\})$. Let 
$T_1:C(X) \rightarrow E$ be the isometric embedding defined by $T_1(f)(x,a)=f(x)$ and $T_1(f)(x,b)=0$.

We now consider a linear isometry extension operator $T_2:E \rightarrow C(X \times Y)$ (with $T_2(1)=1$)
given by Borsuk's Theorem, see \cite{FHHMZ} p. 250 and \cite{B}. Set $T=T_2 \circ T_1$; then $T$ is the desired operator.
Indeed, it is clear that $T$ is an isometry so that $T(f) \diagup X \times \{a\}=T_1(f)$. Let $g \in E:
g \diagup X \times \{a\}=0$ and $||g||=1$; then it is easy to see that the set $T(S) \cup \{T_2(g)\}$ is a 1-\eq
subset of $C(X \times Y)$.
\end{proof}

\noindent \textbf{Remark 4} A compact Hausdorff space $X$ is said to be a Dugundji space, if for any compact
space $Y$ with $X \subseteq Y$ there is a linear extension operator $T:C(X) \rightarrow C(Y)$; that is, a linear 
operator T such that:\\
(i) $T(1)=1$, (ii) $||T||=1$ and (iii) $T(f) \diagup X=f$, for $f \in C(X)$ (see \cite{B}).\\
It is clear that such an operator is an isometry. We note that: (a) Every compact metric space is Dugundji and
(b) the class of Dugundji spaces is closed under cartesian products.

It follows in particular from the above remark that Lemma 5 remains true if we assume that X is any 
Dugundji space and $Y$ any compact space with at least two points.\\

In order to prove Theorem 4 we need to introduce some notation and to remind the reader of some concepts.
Let $\{X_\gamma:\gamma \in \Gamma\}$ be a family of topological spaces. Set $X=\prod\limits_{\gamma \in \Gamma} X_\gamma$,
endowed with the product topology. Fix a point, say $O=(o)_{\gamma \in \Gamma}$ in $X$. For any $\emptyset \neq
A \subseteq \Gamma$ we set 
$ X_A=\prod\limits_{\gamma \in A} X_\gamma \times \{o\}_{\gamma \in \Gamma \setminus A}$
and define a map $\pi_A:X \rightarrow X$ by $\pi_A(x)(\gamma) = \left\{ \begin{array}{ll}
x(\gamma) & ,\gamma \in A\\ \, o & ,\gamma \in \Gamma \setminus A\\
\end{array} \right.$. The map $\pi_A$ is a continuous retraction with $\pi_A(X)=X_A$, which in its
turn induces a norm one projection $P_A:C(X) \rightarrow C(X)$ by the rule, $P_A(f)(x)=f(\pi_A(x))$,
for $f \in C(X)$ and $x \in X$. Note that the range of $P_A$ is identified with the range of the isometry
$T_A:C(X_A) \rightarrow C(X)$ defined by $T_A(g)=g \circ \pi_A$, for $g \in C(X_A)$.

We remind the reader that a map $f:Y \subseteq X \rightarrow \Bbb{R}$ depends on a set $A \subseteq \Gamma$, if
whenever $x,y \in Y$ and $\pi_A(x)=\pi_A(y)$ then $f(x)=f(y)$. It is well known that if each $X_\gamma$ is separable,
then every continuous map $f:X \rightarrow \Bbb{R}$ depends on countably many coordinates. In fact there is a countable
set $A \subseteq \Gamma$ and a continuous map $h:X_A \rightarrow \Bbb{R}$ such that $f=h \circ \pi_A$ (see \cite{E} pp.157-9).

\beg \textbf{(of Theorem 4)}\\
We may (and shall) assume that $|\Gamma| \ge \omega_1$ (if $|\Gamma| \le \omega$, then $X$ is compact metric
and the result holds true). So let $S$ be any \eq set in $C(X)$ with $|S|<|\Gamma|$. Since by Theorem 3,
$m(C(K)) \ge \omega$ for any infinite compact space $K$, we may also assume that $S$ is infinite. We are going
to prove that $S$ can be extended (cf. the proof of Prop. 2). Since each continuous function on $X$ depends on
countably many coordinates and $|S|<|\Gamma|$, it follows that there is an $A \subseteq \Gamma$ with
$|A|=|S|<|\Gamma|$ such that each member of $S$ depends on A. So if $A \subseteq B \subseteq \Gamma$, then
the set $P_B(S)(=S)$ is an \eq subset of $C(X_B)$. Let $\gamma_0 \in \Gamma \setminus A$; set $B=A \cup \{\gamma_0\}$.
Then $P_A(S)$ is an \eq set in $C(X_A)$, $X_A$ is a Dugundji space and $X_B=X_A \times X_{\gamma_0}$, therefore 
by Lemma 5 and Remark 4 we can extend $S$ to an \eq subset $S \cup \{f\}$ of $C(X_B) \subseteq C(X)$. The proof of the theorem
is complete.
\end{proof}

Let $T$ be the dyadic tree, i.e. $T=\bigcup\limits_{n=0}^\infty \{0,1\}^n$ ordered by the relation "$s$ is an initial segment of $t$",
denoted by $s \le t$. By the term chain (resp. antichain) of $T$ we mean a set of pairwise comparable (resp. incomparable)
elements of $T$. A branch of $T$ is any maximal chain of $T$.

Let $A=\{s_1,s_2,\ldots,s_n,\ldots\}$ be any infinite antichain of $T$. We let for any $n \in \Bbb{N}$,
$A_n=\{k \le |s_n|:s_n(k)=1\}$ and $B_n=\{k \le |s_n|:s_n(k)=0\}$  (where $|s|$ denotes the length of $s \in T$,
that is, if $s=(s(1),s(2),\ldots,s(m))$, then $|s|=m$).

\begin{lemm}
The family $\mathcal{F}(A)=\{(A_n,B_n):n \in \Bbb{N}\}$ is a linked family of (finite) pairs of $\Bbb{N}$.
\end{lemm}

\beg
Let $n,m \in \Bbb{N}$ with $n<m$; we then have that $s_n$ and $s_m$ are incomparable. Let $k \le min \{|s_n|,|s_m|\}$
such that $s_n(k) \neq s_m(k)$. If $s_n(k)=1$ then $s_m(k)=0$, thus $k \in A_n \cap B_m$.
If $s_n(k)=0$ then $s_m(k)=1$ and thus $k \in A_m \cap B_n$. So we are done.
\end{proof}

\noindent \textbf{Example 2}  Let $\Delta=\{0,1\}^\Bbb{N}$ be the Cantor set. If $\sigma=(\sigma_1,\sigma_2,\ldots,
\sigma_n,\ldots) \in \Delta$, then the sequence \[C(\sigma)=\{(\sigma_1),(\sigma_1,\sigma_2),\ldots,
(\sigma_1,\sigma_2,\ldots,\sigma_n),\ldots\}\] is a maximal chain (that is, a branch) of $T$ and the sequence
 \[A(\sigma)=\{(\varphi(\sigma_1)),(\sigma_1,\varphi(\sigma_2)),\ldots,(\sigma_1,\ldots,\sigma_{n-1},\varphi(\sigma_n)),\ldots\},\]
where $\varphi(0)=1$ and $\varphi(1)=0$, is a maximal antichain of $T$.

It follows immediately from Lemma 6 that the family $\mathcal{F}(A(\sigma))$ is a linked family of pairs of $\Bbb{N}$.
Furthermore, the family $\mathcal{F}(A(\sigma)) \cup \{(A_\omega,B_\omega)\}$, where $A_\omega=\{n \in \Bbb{N}:\sigma_n=1\}$
and $B_\omega=\{n \in \Bbb{N}:\sigma_n=0\}$ is linked. Indeed, let $(t_n)$ be an arbitrary sequence in the interval $(0,1)$
and let \[p_1=(\varphi(\sigma_1),t_1,\ldots,t_1,\ldots)\]
\[p_2=(\sigma_1,\varphi(\sigma_2),t_2,\ldots,t_2,\ldots)\]
\[\vdots\]
\[p_n=(\sigma_1,\ldots,\sigma_{n-1},\varphi(\sigma_n),t_n,\ldots,t_n,\ldots)\]
\[\vdots\]
\[p_\omega=(\sigma_1,\ldots,\sigma_n,\ldots)=\sigma.\]
Set $S=\{p_n:n \ge 1\}$, then it is obvious that $S \subseteq [0,1]^\Bbb{N} \cap c$, where $c$ is the \ban
of real convergent sequences and that $S \cup \{p_\omega\} \subseteq [0,1]^\Bbb{N} \cap \ell_\infty$. It is also
clear that $S \cup \{p_\omega\}$ is a 1-\eq subset of $\ell_\infty \cong C(\beta \Bbb{N})$ and that the linked
family of pairs of $\Bbb{N}$ corresponding to $S \cup \{p_\omega\}$ according to Lemma 1 is the family
$\mathcal{F}(A(\sigma)) \cup \{(A_\omega,B_\omega)\}$ (cf. also Remarks 1 (2) and (3) and recall that the \ban
$\ell_\infty$ is isometric to $C(\beta \Bbb{N})$, where $\beta \Bbb{N}$ is the Stone-$\check{C}$ech compactification
of the discrete set $\Bbb{N}$).

\begin{lemm}
For any $\sigma \in \Delta$ the set $S \cup \{p_\omega\}$ is a maximal \eq subset of $\ell_\infty$, More exactly:

(a) If the sequence $\sigma$ is eventually constant, then $S \cup \{p_\omega\}$ is a maximal \eq subset of $c$ and

(b) If $\sigma$ is not eventually constant, then $S$ is a maximal \eq subset of $c$.
\end{lemm}

\beg Let $\alpha=(\alpha_n) \in \ell_\infty$ such that $||\alpha-p_n||_\infty=1$, for all $n \ge 1$.
From the equations $||\alpha-p_1||_\infty=||\alpha-p_2||_\infty=1$, we get that $|\alpha_1-\varphi(\sigma_1)| \le 1$
and $|\alpha_1-\sigma_1| \le 1$; since $\{\sigma_1,\varphi(\sigma_1)\}=\{0,1\}$ we conclude that $\alpha_1 \in [0,1]$.
From the equations $||\alpha-p_2||_\infty=||\alpha-p_3||_\infty=1$, we get that $|\alpha_2-\varphi(\sigma_2)| \le 1$
and $|\alpha_2-\sigma_2| \le 1$; since $\{\sigma_2,\varphi(\sigma_2)\}=\{0,1\}$ we conclude that $\alpha_2 \in [0,1]$.
Similarly we conclude that $\alpha_n \in [0,1]$, for all $n \ge 1$.

Since we have $|t_1-\alpha_k|<1$, for all $k \ge 2 $, we get that $1=||p_1-\alpha||_\infty=|\alpha_1-\varphi(\sigma_1)|$,
which clearly implies that $\alpha_1=\sigma_1$. We observe that $|t_2-\alpha_k|<1$, for all $k \ge 3$, therefore
$1=||p_2-\alpha||_\infty=max\{|\alpha_1-\sigma_1|,\,|\alpha_2-\varphi(\sigma_2)|\}=|\alpha_2-\varphi(\sigma_2)|$, which
implies that $\alpha_2=\sigma_2$. In the same way we get that $\alpha_n=\sigma_n$, for all $n \in \Bbb{N}$. So we are done.
\end{proof}

\begin{cor}
The family $\mathcal{F}(A(\sigma)) \cup \{(A_\omega,B_\omega)\}$ is a maximal linked family of pairs of $\Bbb{N}$,
for any $\sigma \in \Delta$.
\end{cor}

\beg It follows immediately from Lemma 7 and Remarks 1 (2) and (3).
\end{proof}

\begin{theo}
$m(\ell_\infty)=\omega(=m(c))$.
\end{theo}

\beg It is an immediate consequence of Theorem 3 and Lemma 7 (The fact that $m(c)=\omega$ also follows
from Theorem 3, this is so because $c$ is isometric to the space $C(\widetilde{\Bbb{N}})$, where $\widetilde{\Bbb{N}}$
is the one point compactification of the discrete set $\Bbb{N}$.)
\end{proof}

\begin{theo}
Let $K$ be a compactification of the discrete set $\Bbb{N}$. Then we have that $m(C(K))=\omega$.
\end{theo}

\beg Let $\sigma=(\sigma_1,\sigma_2,\ldots,\sigma_n,\ldots) \in \Delta$ be an eventually constant sequence,
that is, there are $N \in \Bbb{N}$ and $i \in \{0,1\}$ such that $\sigma_{N+\lambda}=i$, for $\lambda \ge 1$.
Also let $(t_n) \subseteq (0,1)$ be an arbitrary sequence. We define a 1-\eq subset of the space $C(K)$ as follows:
\[f_1=\varphi(\sigma_1) \cdot \chi_{\{1\}}+t_1 \cdot \chi_{K \setminus \{1\}}\]
\[f_2=\sigma_1 \cdot \chi_{\{1\}}+\varphi(\sigma_2) \cdot \chi_{\{2\}}+t_2 \cdot \chi_{K \setminus \{1,2\}}\]
\[\vdots\]
\[f_n=\sigma_1 \cdot \chi_{\{1\}}+\ldots+\sigma_{n-1} \cdot \chi_{\{n-1\}}+\varphi(\sigma_n) \cdot \chi_{\{n\}}+
t_n \cdot \chi_{K \setminus \{1,2,\ldots,n\}}\]
\[\vdots\]
\[f_\omega=\sigma_1 \cdot \chi_{\{1\}}+\ldots+\sigma_N \cdot \chi_{\{N\}}+i \cdot \chi_{K \setminus \{1,2,\ldots,N\}}\]
Let $\pi:\beta \Bbb{N} \rightarrow K$ be the continuous surjective map so that $\pi(n)=n$, for $n \in \Bbb{N}$. Then the
operator $\varPhi:C(K) \rightarrow C(\beta \Bbb{N}) \cong \ell_\infty$ defined by $\varPhi(f)=f \circ \pi$, for $f \in C(K)$,
is an isometric embedding. Since clearly $\varPhi(f_n)=p_n$, for $n \le \omega$, where $\{p_n: n \le \omega\}$ is the sequence
defined in Lemma 7, we get the conclusion.
\end{proof}

The aforementioned results (Theorems 5 and 6) culminate in the following more general result.

\begin{theo}
Let $\Gamma$ be an infinite set and $K$ be a compact space containing an infinite set of isolated points 
having a unique limit point. Then $m(\ell_\infty(\Gamma))=m(C(K))=\omega$.
\end{theo}

\beg Let $D=\{x_n:n \in \Bbb{N}\}$ be a sequence of distinct isolated points of $K$. We consider a non
eventually constant sequence $\sigma=(\sigma_1,\ldots,\sigma_n,\ldots) \in \Delta$ and a dense sequence
$(t_n) \subseteq (0,1)$. We define a 1-\eq subset of $\ell_\infty(K)$, as follows:
\[f_n(x)= \left\{ \begin{array}{ll} \sigma_k & ,x=x_k \, ,k \le n-1 \\ \varphi(\sigma_n) & ,x=x_n \\ 
t_n &,x \in K \setminus \{x_1,x_2,\ldots,x_n\}\\ \end{array} \right.\,\,\,\,
f_\omega(x)= \left\{ \begin{array}{ll} \sigma_n &,x=x_n \, ,n \in \Bbb{N} \\ t_1 & ,x \in K \setminus D\\
\end{array}. \right.\]
It is clear that the set $S=\{f_n:n \in \Bbb{N} \} \subseteq C(K)$.

\begin{claim}
The set $S$ is a maximal \eq subset of $C(K)$ and the set $S \cup \{f_\omega\}$ is a maximal \eq
subset of $\ell_\infty(K)$. 
\end{claim}

\underline{Proof of the Claim:} Let $f \in \ell_\infty(K)$ such that 
\begin{equation}
||f-f_n||_\infty=1,\, \textrm{for} \,\,n \in \Bbb{N}.
\end{equation}
It then follows from (1) that $|f(x)-f_n(x)|=|f(x)-t_n| \le 1$ for all $x \in K \setminus D$ and $n \in \Bbb{N}$. Since $t_n$
is dense in $(0,1)$ we get that $0 \le f(x) \le 1$, for $x \in K \setminus D$. We also have from (1) that $|f(x_n)-\sigma_n| \le 1$
and $|f(x_n)-\varphi(\sigma_n)| \le 1$, for $n \in \Bbb{N}$; so we get that $0 \le f(x_n) \le 1$, for $n \in \Bbb{N}$.
Therefore, 
\begin{equation}
0 \le f(x) \le 1,\, \textrm{for} \, x \in K.
\end{equation}

Let now $N \in \Bbb{N}$; we then have $f_N(x)=t_N$, for $x \in K \setminus D$. It follows from (2) that
\[|f(x)-f_N(x)|=|f(x)-t_N| \le max\{t_N,1-t_N\}<1 \, \textrm{for} \, x \in K \setminus D. \]
So we get that \[||f-f_n||_\infty=\sup_{x \in D} |f(x)-f_n(x)|=1, \, \textrm{for} \, n \in \Bbb{N}.\]
The method of proof of Lemma 7 then yields that $f(x_n)=\sigma_n$, for $n \in \Bbb{N}$.

Assume for the moment that $D$ has a unique limit point, say $x$, thus $x_n \rightarrow x$. Since $f(x_n)=\sigma_n=1$
for infinite $n \in \Bbb{N}$ and $f(x_n)=\sigma_n=0$ for infinite $n \in \Bbb{N}$, we conclude that $f$ cannot be
continuous on $K$. So the set $S$ is a maximal \eq subset of $C(K)$.

Also $S \cup \{f_\omega\}$ is a maximal \eq subset of $\ell_\infty(K)$, since (for an arbitrary infinite set $K$)
\[||f-f_\omega||_\infty=\sup_{x \in K \setminus D} |f(x)-f_\omega(x)|=\sup_{x \in K \setminus D} |f(x)-t_1|
\le max\{t_1,1-t_1\}<1.\]
The proof of the Claim and hence of the Theorem is complete.
\end{proof}

\begin{cor}
Let $K$ be an infinite compact scattered space. Then $m(C(K))=\omega$.
\end{cor}

\beg Since $K$ is scattered, the set of isolated points of $K$ is dense in $K$; moreover $K$ is sequentially
compact. So let $(x_n)$ be a sequence of distinct isolated points of $K$, then $(x_n)$ has a convergent
subsequence, say $y_n=x_{m_n}$, $n \in \Bbb{N}$. Therefore the set $D=\{y_n:n \in \Bbb{N}\}$ has a unique
limit point and the previous theorem can be applied.
\end{proof}

\noindent \textbf{Remarks 5} (1) Let $X,\,Y$ be compact spaces, $\pi:X \rightarrow Y$ a continuous surjective
map which is non irreducible (i.e., there is $\Omega \subsetneqq X$ compact such that $\pi(\Omega)=Y$) and 
$T:C(Y) \rightarrow C(X)$ be the linear isometry given by $T(f)=f \circ \pi$, for $f \in C(Y)$. We consider
a 1-\eq subset $S$ of $S^{+}_{C(Y)}$, then it is rather easy to prove that there is $g \in S^{+}_{C(X)}$,
such that the set $T(S) \cup \{g\}$ is equilateral.

Given this result, it can be shown by transfinite induction, that if a compact non metrizable space $K$ is roughly
the "limit" of  a long system of "smaller" compact spaces connected by non irreducible maps, then $m(S^{+}_{C(K)} \cup \{0\},0)
\ge \omega_1$. This is the case for instance when: 

(a) $K$ is an Eberlein, or more general a Corson compact space, because then $K$ admits a long sequence of compatible
retractions, see \cite{N}, \cite{DGZ} and

(b) $K$ is a compact group, since then $K$ is a projective  limit of compact metrizable groups, see \cite{HR}.\\
(2) Note that, given any infinite set $\Gamma$, the Banach spaces $c_0(\Gamma)$ and $c(\widetilde{\Gamma})$ are
isomorphic, where $\widetilde{\Gamma}$ is the one point compactification of the discrete set $\Gamma$. But if $\Gamma$
is uncountable, then by Prop. 2, Cor. 4 and the above remark we have
\[m(C(\widetilde{\Gamma}))=\omega<m(S^{+}_{C(\widetilde{\Gamma})} \cup \{0\},0)=m(c_0(\Gamma))=|\Gamma|.\]
\flushleft (3) The following example is related with Remark 1(3):\\
Let $K$ be a compact nonempty space. We denote by $\Omega$ the disjoint union of the compact spaces
$K$ and $\widetilde{\Bbb{N}}$, where $\widetilde{\Bbb{N}}=\Bbb{N} \cup \{\infty\}$ is the one point compactification of 
$\Bbb{N}$. Let $\sigma=(\sigma_1,\sigma_2,\ldots,\sigma_n,\ldots) \in \Delta$ be an eventually constant sequence, so that
$\sigma_n=i \in \{0,1\}$ for $n \ge N$. We define a 1-\eq set $S=\{f_n:n \le \omega\} \subseteq [0,1]^\Omega \cap C(\Omega)$
as follows:
\[f_n(x)= \left\{ \begin{array}{ll} \sigma_k & ,x=k \, ,k \le n-1 \\ \varphi(\sigma_n) & ,x=n \\ 
\frac{1}{2} &,x \in \Omega \setminus \{1,2,\ldots,n\}\\ \end{array} \right. ,n \in \Bbb{N}\,\,\textrm{and}\,
f_\omega(x)= \left\{ \begin{array}{ll} \sigma_n &,x=n \\ i & ,x=\infty \\ \frac{1}{2} & ,x \in K \\
\end{array}. \right.\]
It is easy to verify that the linked family $\mathcal{F}$ corresponding to this \eq set according to Lemma 1 is maximal
(cf. Lemma 7) and also that $S$ can be extended to an \eq set in $C(\Omega)$.

\noindent S.K.Mercourakis, G.Vassiliadis\\
University of Athens\\
Department of Mathematics\\
15784 Athens, Greece\\
e-mail: smercour@math.uoa.gr

\hspace{0.7cm} georgevassil@hotmail.com
\end{document}